\documentclass[12pt]{article} 
\usepackage{graphicx,amsmath, amssymb, verbatim}
\usepackage[all]{xy}
\newtheorem{proposition}{Proposition}

\newcommand {\C } {\mathbb{C}} 
\newcommand {\Q } {\mathbb{Q}} 

\newcommand {\p } {\mathbb{P}} 
\newcommand {\R } {\mathbb{R}} 
\newcommand {\Z} {\mathbb{Z}}
\newcommand {\oo} {\mathcal O}

\newcommand{\ba}{\begin{eqnarray}}
\newcommand{\ea}{\end{eqnarray}}
\newcommand{\bmm}{\begin{pmatrix}}
\newcommand{\emm}{\end{pmatrix}}
\newcommand{\no}{\nonumber}

\begin{document} 

\title{\bf{On equivariant mirror symmetry for local $\p^2$}}

\author{Brian Forbes \\
\it Research Institute for Mathematical Sciences \\ \it  Kyoto University \\ \it  Kyoto 606-8502, Japan \\ \it{brian@kurims.kyoto-u.ac.jp}\\ \\
Masao Jinzenji  \\ \it Division of 
Mathematics, Graduate School of Science \\ \it Hokkaido University \\ 
\it  Sapporo  060-0810, Japan\\
\it jin@math.sci.hokudai.ac.jp } 

\maketitle

\begin{abstract}
We solve the problem of equivariant mirror symmetry for
 $K_{\p^2}=\oo(-3)\rightarrow \p^2$ for the (three) cases of one independent
 equivariant parameter. This gives a decomposition of mirror symmetry for
 $K_{\p^2}$ into that of three subspaces, each of which may be considered
 independently. Finally, we give a new interpretation of
 mirror symmetry for $\oo(k)\oplus \oo(-2-k)\rightarrow \p^1$.  
\end{abstract}

\section{Introduction}

\subsection{Background}

In recent years, our mathematical understanding of mirror symmetry has
progressed dramatically (see in particular \cite{CG,J,I,FJ3,FJ4} for
results relevant to the present paper). Nonetheless,
there remains an entire theory of mirror symmetry which does not
evidently fall under this umbrella, and which remains undeveloped;
we shall call this \it equivariant mirror symmetry. \normalfont

Of course, as any theory of mirror symmetry must, our tentative `equivariant
mirror symmetry' has its roots in equivariant Gromov-Witten
theory. What this means is essentially that, in addition to the usual
data coming from the Gromov-Witten calculation, one must also keep track of
the torus weights corresponding to the natural torus action on our
(toric) Calabi-Yau space.
The first example of this was detailed in \cite{BP}, where they solved
equivariant Gromov-Witten theory for arbitrary rank 2 bundles over a
genus $g$ Riemann surface $\Sigma_g$, where one allows arbitrary torus
weights $(\mu_1,\mu_2)$ to act on the bundle.

Using this example as a starting point, for the special case
$\Sigma_0=\p^1$, we sucessfully worked out the equivariant mirror
computation in \cite{FJ3,FJ4}.
Namely, it was shown in \cite{FJ3} that by appropriately
adapting the results of \cite{CG}, the genus zero invariants of
\cite{BP} could be computed via mirror symmetry. This approach was
simplified and streamlined in \cite{FJ4}, where a formula for the mirror
map agreeing with the physical analysis of \cite{M} was derived (as well
as Picard-Fuchs equations and Yukawa couplings).

The above then naturally leads us to the question of how to describe
equivariant mirror symmetry for the next most complicated example:
$K_{\p^2}=\oo(-3)\rightarrow \p^2$. Although the only difference between this and
the examples of $\cite{FJ3,FJ4}$ is that now $b_4(K_{\p^2})\ne 0$, as we will
see below, even this small difference leads to a vastly more challenging
problem. Nonetheless, we find that an approach which is more or less the
natural generalization of \cite{FJ3,FJ4} allows us to deal with the
mirror symmetry problem for this space, in the case that we
have only one independent equivariant parameter. We check our results via the (Chern-Simons)
computations of \cite{AMV} and find agreement with their results. As
the Chern-Simons theory contains significantly more data than ordinary
mirror symmetry, our calculations represent a substantial
refinement of mirror symmetry. Remarkably,
we also
find a method of modifying the usual localization calculation for
$K_{\p^2}$ which exactly reproduces the `refined Gromov-Witten
invariants' of \cite{AMV}.

In addition, we find a new and simplified approach to deriving the
results of \cite{FJ4}, which allows us to produce closed formulas for
the two-curve formula (defined below). 

The organization of this paper is as follows. Section 2 briefly outlines
the approach of \cite{AMV} to the equivariant mirror problem. Section 3
establishes the equivariant mirror calculation for $K_{\p^2}$ for one
independent equivariant parameter (this gives 3 separate cases; see
Figure 1) both by directly using \cite{FJ3}, and by applying the
simplifications of \cite{FJ4}. Section 4 gives a novel approach to
localization which allows us to directly work out the fully equivariant
Gromov-Witten theory of $K_{\p^2}$. Finally, Section 5 presents our new
approach to the problem of $\oo(k)\oplus \oo(-2-k)\rightarrow
\p^1$. Relevant facts about $I$ and $J$ functions and their relation to
mirror symmetry is presented in the Appendix.

\subsection{Summary of results}

\begin{figure}[t]
\label{all3}
\centering
\input{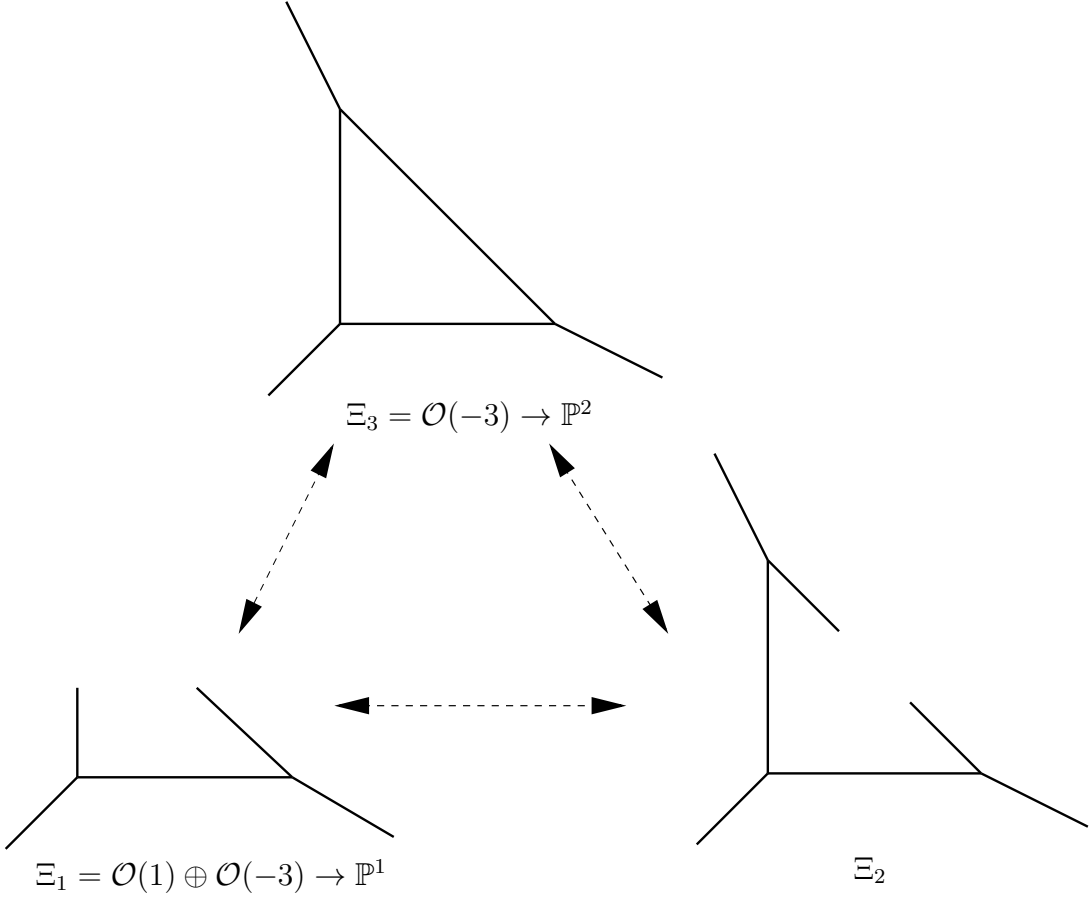}
\caption{Decomposition of $K_{\p^2}$ via equivariant parameters}
\end{figure}

In this paper, we work out equivariant mirror symmetry for
$\oo(-3)\rightarrow \p^2,$ for the case of one independent equivariant parameter. The geometric meaning of this is depicted in
Figure 1:  certain choices of equivariant  weights decompose $\oo(-3)\rightarrow \p^2$ into
a family of three subspaces $\Xi_1=\oo(1)\oplus \oo(-3)\rightarrow
\p^1,\Xi_2$ and $\Xi_3=\oo(-3)\rightarrow \p^2$, and we give mirror
symmetry constructions which reproduce the Gromov-Witten invariants of each.

Concretely, we find the following:
\begin{proposition}
Mirror symmetry for the spaces $\Xi_1,\Xi_2$ and $\Xi_3$ is the same as
 that of
\ba
\Xi_1&:& \ \ \ 
\oo(1)_{\nu}\oplus \oo(-3)_{-\nu} \rightarrow \p^1_{[0,0]} \\ \no  \\ \no
\Xi_2&:& \ \ \ 
\oo(1)_{\nu}\oplus \oo(-3)_{-2\nu} \rightarrow \p^1_{[0,\nu]} \\ \no  \\ \no
\Xi_3&:& \ \ \
 \oo(-3)_{-3\nu} \rightarrow \p^2_{[\nu,\nu,\nu]}
\ea
where the subscripts denote equivariant weights acting on each
 term. Equivalently, this can be represented by the following toric data
 subject to the cohomology relations: 
\ba
 \begin{pmatrix}
1&1&1&-3\\
&&\nu&-\nu
\end{pmatrix} \ \ \ &:& p^2=0
\\ \no
\begin{pmatrix}
1&1&1&-3\\
&\nu&\nu&-2\nu
\end{pmatrix} \ \ \ &:& p(p+\nu)=0
\\ \no
\begin{pmatrix}
1&1&1&-3\\
\nu&\nu&\nu&-3\nu
\end{pmatrix} \ \ \ &:& (p+\nu)^3=0
\ea
Above, the second row of
 each matrix specifies the equivariant weights, and $\begin{pmatrix}
1&1&1&-3
\end{pmatrix}$ is the standard toric vector for $\oo(-3)\rightarrow \p^2$.
\end{proposition}

We use this to confirm a prediction of \cite{AMV} regarding the relationship of
the so-called `refined Gromov-Witten invariants' computed via
Chern-Simons theory in that paper, and equivariant Gromov-Witten
invariants. 

We expect that the eventual `full' equivariant mirror symmetry
computation (which will give a Gromov-Witten generating function
agreeing with \cite{AMV}) will use our methods to some extent. However,
at present, the complicated interplay of Birkhoff factorized and
non-Birkhoff factorized functions involved have put the full theory out
of reach.

The second result of this paper is a new interpretation of the mirror
symmetry calculation of $\oo(k)\oplus \oo(-2-k)\rightarrow \p^1$ with
antidiagonal equivariant weights $(\nu,-\nu)$ acting on the bundle. We
find
\begin{proposition}
Mirror symmetry for $X_k=\oo(k)\oplus \oo(-2-k)\rightarrow \p^1$ with the
 antidiagonal action (the equivariantly Calabi-Yau case) can be
 determined from the toric data 
\ba
\begin{pmatrix}
1&1&1&-1&-1&-1\\
&&\nu/k&-\nu/k&-\nu&-\nu
\end{pmatrix}
\ea
which corresponds to the geometry
\ba
\oo(1)_{\nu/k}\oplus\oo(-1)_{-\nu/k}\oplus\oo(-1)_{-\nu}\oplus\oo(-1)_{-\nu}\rightarrow \p^1
\ea
where the subscripts denote equivariant weights (and we are considering
 the trivial action on $\p^1$). This reproduces the mirror maps of \cite{FJ4}\cite{M} and gives
 a closed formula for the generating function of Gromov-Witten
 invariants of $X_k$:
\ba
W_k(x)=\sum_{j=1}^{\infty}\frac{(-1)^{(j-1)k}}{j^2(j-1)!}\prod_{m=1}^{j-1}(k(2+k)j+m)x^k
\ea
\end{proposition}

Finally, we give a generalization of the above result which gives closed
formulas for the mirror map and GW invariants of a large class of geometries.

\bigskip

\bigskip

\bf{Acknowledgements}\normalfont

We owe special thanks to Marcos Marino for suggesting this problem, and
for helpful correspondence. We would also like to thank C. Doran,
S. Hosono, H. Iritani, B. Kim and 
K. Saito for helpful discussions.

\section{The idea of equivariant mirror symmetry}

Here, we would like to outline the general principles that guide our
study of equivariant mirror symmetry. 

We recall Section 7.6 of \cite{AMV} and its application to our current
problem, $K_{\p^2}$. The idea there was to define a set of `refined'
Gromov-Witten invariants of $K_{\p^2}$ by setting the K\"ahler parameter
of each $\p^1 \hookrightarrow \p^2$ to a different value, which is
apparently a natural operation from the vantage of the dual Chern-Simons
theory. They then went
on to compute these new Gromov-Witten invariants as rational functions
of the
three K\"ahler parameters. The expectation is that these should
somehow correspond to the equivariant Gromov-Witten theory of
$K_{\p^2}$, with equivariant weights given as \cite{M2}:
\ba
\label{fullyequiv}
\begin{pmatrix}
1&1&1&-3 \\
\mu_1&\mu_2&\mu_3&-\mu_1-\mu_2-\mu_3
\end{pmatrix}.
\ea 
While these two concepts are not readily seen to be the same, in this
paper, for the cases $(\mu_1,\mu_2,\mu_3)=(0,0,\nu),(0,\nu,\nu)$ and $(\nu,\nu,\nu)$,
we find agreement with the results of \cite{AMV}.

A general formulation of what is being computed in \cite{AMV} would go something as
follows. Let $X$ be a noncompact toric Calabi-Yau threefold, and let $\Gamma$ be the
image of the moment map, i.e. the toric graph associated to $X$. We denote the set of bounded edges of $\Gamma$
by $\{\gamma_i \}$, and
associate a parameter $\nu_i$ to each $\gamma_i$. Then the generating
function of Gromov-Witten invariants of \cite{AMV} will be a function
\ba
W(\nu_1 \dots \nu_n) \in \Q[\nu_1 \dots  \nu_n][[x]]
\ea
such that e.g.
$W(\nu_1,
\nu_2=0 \dots \nu_n=0)$ is the generating function of Gromov-Witten
invariants of maps 
\ba
f:\p^1\rightarrow X, \ \ \ f_*(\p^1)\in \gamma_1.
\ea
Similarly, $W(\nu_1,
\nu_2, \nu_3=0 \dots \nu_n=0)$ will count maps from $\p^1$ into
$\gamma_1 \cup \gamma_2$, etc. Then mirror symmetry in this
sense would mean that there would be a mirror manifold $Y_{\nu_1 \dots
\nu_n}$ to $X$ such that the period integrals of $Y$ would allow us to
compute $W(\nu_1 \dots \nu_n)$. According to this picture, in addition to ordinary
mirror symmetry, equivariant mirror symmetry should also include mirror
symmetry for each 
\ba
\big(N_{C/X}\rightarrow C\big) \hookrightarrow X
\ea
for each $C\hookrightarrow X$,
and for each pair of curves in $X$, etc. Moreover, all of
these should be incorporated into a single picture with the use of the
parameters $\nu_i$.

\section{Equivariant mirror symmetry for $K_{\p^2}$:  one independent equivariant
 parameter}

We now present the techniques for the equivariant mirror
symmetry calculation corresponding to Figure 1. The methods are
surprisingly simple, yet require some highly non-trivial observations in order
to be successfully carried out. 

\subsection{Review of $\Xi_1=\oo(1)\oplus \oo(-3)\rightarrow \p^1$}

We first review the content of \cite{FJ3}, which established the
mirror symmetry calculation for $\Xi_1=\oo(1)\oplus \oo(-3)\rightarrow \p^1$. Recall that the toric data
describing $K_{\p^2}$ is just
\ba
\label{p2vect}
\begin{pmatrix}
1&1&1&-3
\end{pmatrix},
\ea
which is a shorthand notation for the toric quotient
\ba
K_{\p^2}=\{(z_i)\in \C^4 : |z_1|^2+|z_2|^2+|z_3|^2-3|z_4|^2=r\}/S^1
\ea
where $r \in \R^+$ and 
\ba
S^1:(z_1 \dots z_4)\longrightarrow (e^{i \theta}z_1,e^{i \theta}z_2,e^{i
\theta}z_3,e^{-3 i \theta}z_4).
\ea

From this, we see immediately that the three 1's in (\ref{p2vect}) represent homogeneous
coordinates on $\p^2$. So, in order to get at the geometry of $\oo(1)\oplus
\oo(-3)\rightarrow \p^1$, we want to treat $z_3$ as a fiber variable
with respect to $\p^1=[z_1,z_2]$. This can be accomplished by
considering a partially equivariant theory on $K_{\p^2}$:
\ba
\label{pep2vect}
\begin{pmatrix}
1&1&1&-3 \\
&&\nu&-\nu
\end{pmatrix}
\ea
where $\nu \in \C^*$ and the second column represents the torus weights
acting on the coordinates of $K_{\p^2}$. As explained below, we have to
subject the $I$ function corresponding to this vector to the cohomology
relation $p^2=0$, which implies that we are considering the geometry
\ba
\oo(1)_{\nu}\oplus
\oo(-3)_{-\nu}\rightarrow \p^1_{[0,0]}.
\ea
Then, by expanding in
$\nu=\infty$ and performing Birkhoff factorization, as advocated in
\cite{CG}, we can indeed realize mirror symmetry for $\oo(1)\oplus
\oo(-3)\rightarrow \p^1$.

Concretely, this goes as follows. The $I$ function corresponding to
(\ref{pep2vect}) is a cohomology-valued hypergeometric series
\ba
I_1=q^{p/\hbar}\sum_{d \ge 0}\frac{\prod_{m=-3d+1}^0(-3p-\nu+m\hbar)}{\prod_{m=1}^d(p+\nu+m\hbar)\prod_{m=1}^d(p+m\hbar)^2}q^d,
\ea
which is annihilated by the equivariant Picard-Fuchs operator
\ba
\label{equivpf0}
\mathcal D_1=\theta^2(\theta+\nu)-q(-3\theta-\nu)(-3\theta-\nu-\hbar)(-3\theta-\nu-2\hbar),
\ \ \ \ \theta=\hbar q\frac{d}{dq}.
\ea
Above, $p$ is the K\"ahler class on $\p^1$, which means that we enforce
the cohomology relation $p^2=0$ in $I$. This is equivalent to
considering just 2 of the 3 solutions of $\mathcal D_1 f=0$:
\ba
I|_{p=0}, \ \ \ \ \frac{d}{dp}I|_{p=0}.
\ea
We then Birkhoff factorize these two to recover a $J$ function, and
subsequently transform by the mirror map, which gives us the instanton
expansion for the $\Xi_1$ geometry:
\ba
\label{onecurvesuper}
W_{\Xi_1}&=&\sum_{k>0}
\frac{(-1)^{k-1}}{k^2(k-1)!}\prod_{j=1}^{k-1}\big(3k+j\big)x^k \\ \no
&=&x-\frac{7}{4}x^2+\frac{55}{9}x^3-\frac{455}{16}x^4+\frac{3876}{25}x^5-\frac{33649}{36}x^6+\dots
\ea
Note that this is the derivative of the prepotential of
$\Xi_1=\oo(1)\oplus \oo(-3)\rightarrow \p^1$, which must be
(logarithmically) integrated once in $x$ to compute integer invariants
via the multiple cover formula.

\subsection{$\Xi_2$: Two $(1,-3)$ curves in $K_{\p^2}$}

\begin{figure}[t]
\label{2curve}
\centering
\input{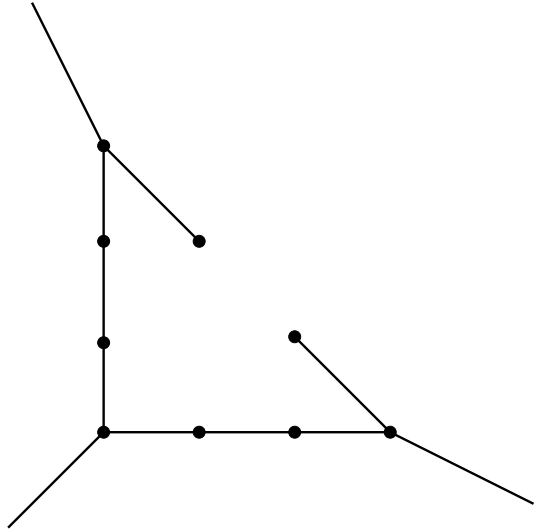}
\caption{The two curve problem $\Xi_2$}
\end{figure}

Before diving into the details, we make some effort at clarifying the meaning of the space
$\Xi_2$, which is represented in Figure 2. The basic intuition is that
we should consider the configuration of two curves, intersecting at
right angles, such that the normal bundle of each is $\oo(1)\oplus
\oo(-3)\rightarrow \p^1$, and also such that the $\oo(-3)$ factor of the
two curves is the same.

Unfortunately, at present, we lack a rigorous mathematical
definition for $\Xi_2$. However, when viewed as a localization problem, the idea behind it is
simple. Ordinary
Gromov-Witten theory of $K_{\p^2}$ computes some number based on the
moduli space of degree $d$ maps
\ba
N_d=\#\{ f:\p^1\rightarrow K_{\p^2}\}.
\ea
Now, in order to define the Gromov-Witten invariants of $\Xi_2$, we just
have to remove certain contributions from this number. Namely, if $C_i, i=1,2,3$
represent the three $\p^1$s in $\p^2$, then we want to remove all curves
which contain a component
mapping to $C_3$.

From the physics point of view, the dual Chern-Simons theory evidently provides some sort of
definition for $\Xi_2$ \cite{AMV}, since \cite{AMV} was able to produce
Gromov-Witten invariants for $\Xi_2$. However, this relies on large
$N$ duality, which has clearly not been proven mathematically. Another
possibility, though equally unproven, would be to define this space by
its equivalence to the theory of open strings with two Lagrangian
submanifolds on $\C^3$ (briefly outlined later in the paper). Note that this is
 equivalent to the two-leg topological vertex, with appropriately chosen framings. 

At any rate, even in absence of a rigorous definition of $\Xi_2$, the space does have known
Gromov-Witten invariants, and hence we can approach the problem by
simply looking for a toric model which reproduces these invariants
through mirror symmetry. We use the above intuitive picture of $\Xi_2$
as a guide in seeking the correct toric data. Remarkably, we find that
the most obvious choice of (equivariant) toric vector precisely
describes mirror symmetry for $\Xi_2$ (though some highly non-trivial
manipulations are required to demonstrate this).

The Gromov-Witten invariants of $\Xi_2$ can be computed from the following
generating function, which is the derivative of the
prepotential for the 2 curve case $\Xi_2$ (note the dramatic
similarity with the $\oo(1)\oplus \oo(-3)\rightarrow \p^1$ prepotential
of eqn (\ref{onecurvesuper})):
\ba
\label{2cp}
W_{\Xi_2}&=&2\sum_{k>0}\frac{(-1)^{k-1}}{k^2(k-1)!}\prod_{j=1}^{k-1}\big(5k+j\big)x_1^k \\
&=&2x_1-\frac{11}{2}x_1^2+\frac{272}{9}x_1^3-\frac{1771}{8}x_1^4+\frac{47502}{25}x_1^5-\frac{162316}{9}x_1^6+\dots
\ea
Note that this is a function of only one K\"ahler parameter. In order to
see a 2 K\"ahler parameter expansion, one would have to use two
different equivariant parameters; however, incorporating this into the
full equivariant mirror calculation of $K_{\p^2}$ is beyond the scope of
the present work.

So, on to the calculation. From the previous section, we saw that the
$\Xi_1$ calculation could be done by use of the toric
data
\ba
\begin{pmatrix}
1&1&1&-3\\
&&\nu&-\nu
\end{pmatrix}
\ea
together with some semi-involved manipulations of the period integrals \cite{CG}\cite{FJ3}.
Thus, the choice which immediately springs to mind for use on $\Xi_2$ is
\ba
\label{2curvevect1}
\begin{pmatrix}
1&1&1&-3\\
&\nu&\nu&-2\nu
\end{pmatrix}
\ea
for mirror symmetry. However, there are major problems with this
that make it seem far too naive to be the right toric data. Namely, on
$\Xi_1=\oo(1)\oplus \oo(-3)\rightarrow \p^1$ above, we could apply the
cohomology relation $p^2=0$ for $\p^1$ to the $I$ function. This
relation then dictates that we must expand the $I$ function in
$\nu=\infty$ \cite{CG}, and we can then use Birkhoff factorization to
exhibit mirror symmetry. However, on $\Xi_2$, since we cannot even
define the space mathematically, we have no guess as to what cohomology
relation to use, and therefore no way of knowing the correct asymptotic
expansion of $\nu$. 

Nonetheless, we will show here that one basic assumption allows us to
essentially use (\ref{2curvevect1}) to exhibit mirror symmetry for
$\Xi_2$. The assumption is, \it  we only need to
use 2 of the 3 solutions of the associated equivariant Picard-Fuchs
system in order to determine mirror symmetry. \normalfont The
justification for this assumption is that, on both $\oo(1)\oplus
\oo(-3)\rightarrow \p^1$ and (equivariant) $\oo(-3)\rightarrow \p^2$, computation of
the prepotential requires only two of three Picard-Fuchs
solutions. Thus, it is reasonable to hope that the `intermediate' case
$\Xi_2$ also enjoys this feature.

We can now see relatively simply that this assumption allows us to derive mirror symmetry for $\Xi_2$. The first implication of our assumption is that we have to
work instead with an equivariant theory given by 
\ba
\label{2curvevect2}
\begin{pmatrix}
1&1&1&-3\\
&\mu&\nu&-2\nu
\end{pmatrix}.
\ea 
The reason for using a new parameter $\mu$ is as follows. The Picard-Fuchs
system associated to (\ref{2curvevect2}) is 
\ba
\label{equivpf0}
\mathcal D_2=\theta(\theta+\mu)(\theta+\nu)+q(3\theta+2\nu)(3\theta+2\nu+\hbar)(3\theta+2\nu+2\hbar).
\ea
Then if we set
\ba
I_2=q^{p/\hbar}\sum_{d \ge 0}\frac{\prod_{m=-3d+1}^0(-3p-2\nu+m\hbar)}{\prod_{m=1}^d(p+\mu+m\hbar)\prod_{m=1}^d(p+\nu+m\hbar)\prod_{m=1}^d(p+m\hbar)}q^d,
\ea
the 3 solutions of $\mathcal D_2 f=0$ are given by 
\ba
I_2|_{p=0}, \ \ \ I_2|_{p=-\mu}, \ \ \ I_2|_{p=-\nu}. 
\ea
Under our assumption, we want to use two of these three solutions for
our mirror symmetry calculation. If we use the first two of these
solutions, this will force us to asymptotically
expand $I_2$ in $\mu=0, \nu=\infty$. This is because,
using only the first two solutions is equivalent to subjecting $I_2$ to the cohomology relation
$p(p+\mu)=0$, which forces a $\mu=0$ expansion. This then implies that the variables involving $\nu$ are `fiber'
variables with respect to the equivariant $\p^1$ given by $p(p+\mu)=0$,
and hence must be expanded in $\nu=\infty$, according to \cite{CG}.

So, in summary, since our assumption requires that $I_2$ be expanded in
$\mu=0,\nu=\infty$, we see the necessity of introducing a new
equivariant parameter $\mu$ compared to our original guess (\ref{2curvevect1}).

We briefly note that the geometry
(\ref{2curvevect2}) subject to the cohomology relation $p(p+\mu)=0$ is
given by
\ba
\oo(1)_{\nu}\oplus \oo(-3)_{-2\nu}\rightarrow \p^1_{[0,\mu]}
\ea

Now that we have fully set up the problem, we have only to carry out the
calculation. We expand $I_2$ as
\ba
\frac{1}{\prod_{m=1}^d(p+\nu+m\hbar)}, \ \  in \ \  \hbar=0, \\
\frac{1}{\prod_{m=1}^d(p+\mu+m\hbar)}, \ \  in \ \  \hbar=\infty.
\ea
By using the equality 
\ba
\mu \frac{d}{dp}I_2|_{p=0}=-I_2|_{p=0}+I_2|_{p=-\mu},
\ea
we can perform Birkhoff factorization of the two solutions $I_2|_{p=0}, \frac{d}{dp}I_2|_{p=0}$. The result is a $J$
function
\ba
J_{\mu,\nu=1}=1+\frac{pt+t_0}{\hbar}+\frac{((3\mu-5)p+2)W}{\hbar^2}+\dots
\ea
where the mirror map $t$ has a form which is unfortunately not defined
in the $\mu \rightarrow 1$ limit  (we have set $\nu=1$ for simplicity):
\ba
\no
t=-\Big(5+6\mu+(\mu-1)^{-1}\Big)q+\Big(-\frac{33}{2}+117\mu-\frac{1}{2}\big[(\mu-1)^{-1}+7(\mu-1)^{-2}+(\mu-1)^{-3}\big]\Big)q^2+\dots
\ea
Remarkably, $W$ does turn out to be well-defined for $\mu=1$:
\ba
\no
W_{\mu,\nu=1}(x)&=&\big(2\mu-4\big)x+\big(26-\frac{53}{2}\mu+6\mu^2\big)x^2-\big(\frac{3028}{9}-\frac{4646}{9}\mu+246\mu^2-36\mu^3\big)x^3\\
&+&\big(\frac{11601}{2}-\frac{95077}{8}\mu+\frac{17363}{2}\mu^2-2664\mu^3+288\mu^4\big)x^4+\dots
\ea
where $x=e^t$.
Then, simply setting $\mu=1$ gives the 2 curve expansion (\ref{2cp})!
\ba
W_{\mu=1,\nu=1}(x)=-2x+\frac{11}{2}x^2-\frac{272}{9}x^3+\frac{1771}{8}x^4+\dots
\ea
Hence, we have accomplished our goal of exhibiting mirror symmetry
 for the space $\Xi_2$.

We also note that we can perform another consistency check of this
calculation, which gives geometric insight into why this method produces
invariants for $\Xi_2$. That is, by setting $\mu=0$ in
$W_{\mu,\nu=1}(x)$, the resulting invariants agree with those of 
\ba
\oo(1)_{\nu}\oplus \oo(-3)_{-2\nu}\rightarrow \p^1_{[0,0]}.
\ea
Of course, this is precisely the $I$ function we get by setting $\mu=0$
into $I_2$, so this is consistent. What this also says is that the $\mu$
parameter is somehow encoding the information of an extra divisor in
$K_{\p^2}$, i.e. without $\mu$, the result collapses into the invariants
of a single
$\p^1\hookrightarrow K_{\p^2}$.  

As a final note, we mention that this computation can clearly be
combined with the $\oo(1)\oplus \oo(-3)\rightarrow \p^1$ calculation of
the previous section (and also with the full $K_{\p^2}$ case, in fact)
simply by using the toric vector
\ba
\begin{pmatrix}
1&1&1&-3\\
&k\mu&\nu&-(k+1)\nu
\end{pmatrix} \ \ : \ \ p(p+k\mu)=0
\ea 
and going through the same steps as above. Then $k=0$ will give $\Xi_1$
and $k=1$ gives $\Xi_2$.

\subsection{$\Xi_3=K_{\p^2}$ with equivariant parameters}

Finally, we give a brief discussion on including $K_{\p^2}$ into the
above picture. Although this is the only of the three spaces for which
standard mirror symmetry actually applies, how to incorporate this into the
above Birkhoff factorization scheme is not immediately clear. This is
because, as we will see, a Birkhoff factorized function cannot be used
to compute $K_{\p^2}$ invariants, which necessitates including \it both
\normalfont the factorized and unfactorized functions together. Nonetheless,
certain coincidences in the Birkhoff factorization make this possible.

Now, in analogy with (\ref{2curvevect1}) above, the natural thing to use is
\ba
\begin{pmatrix}
1&1&1&-3\\
\nu&\nu&\nu&-3\nu
\end{pmatrix},
\ea
and the most naive application of the above machinery would be to try to
work with
\ba
\label{M}
\begin{pmatrix}
1&1&1&-3\\
\mu&\mu&\nu&-3\nu
\end{pmatrix}.
\ea
with $\mu=0,\nu=\infty$. However, we quickly see that this can't be
right, since the computation reveals that the Gromov-Witten invariants
here are the same as those of $\oo(1)_{\nu}\oplus
\oo(-3)_{-3\nu}\rightarrow \p^1_{[0,0]}$. Further investigation shows
that the prepotential $\mathcal F$ of the most general configuration $\oo(1)_{\nu_1} \oplus
\oo(-3)_{\nu_2}\rightarrow \p^1_{[\mu_1,\mu_2]}$ is a function of only
two parameters:
\ba
\label{F}
\mathcal F= \mathcal F \big(\mu_1-\mu_2,\frac{\nu_2}{\nu_1}\big).
\ea
It is not hard to show that the invariants of $K_{\p^2}$ are not
included here. This might have
been guessed at, as the Birkhoff factorization necessarily removes some
information from the fundamental solution of $K_{\p^2}$, and we need all
the data of the original function to compute the right invariants.

Hence, we have to pursue another approach. As we now show, a simple trick actually
allows us to simultaneously express the invariants
corresponding to the $\hbar$ and $1/\hbar$ expansion of the $I$
function. 

For simplicity, we consider only a single parameter in equivariant
cohomology, i.e. $H^*_{\C^*}(\p^2)$. (Note that, by (\ref{F}), this is
equivalent to using (\ref{M}) if we change $-3\nu$ to $-k\nu$, $k \in \Z$.) This is the same setup we had in
Section 3.1:
\ba
\begin{pmatrix}
1&1&1&-3\\
&&\nu&-\nu
\end{pmatrix},
\ea
but the interpretation will be a bit different this time. First, one can
compute readily that, when we expand the $I$ function corresponding to
this vector, $I_1$, in $\hbar=\infty$, we can see the $K_{\p^2}$
Gromov-Witten invariants (even after imposing the $p^2=0$
relation). This is just due to the fact that the $I$ functions for
$K_{\p^2}$ and the $(1,-3)$ curve are the same.  Then, as the $\oo(1)\oplus \oo(-3)\rightarrow \p^1$
invariants were also computed using $p^2=0$, this is simply a matter of
putting both pieces together. 

Ultimately, the connection between the $K_{\p^2}$ calculation and those
of the previous two sections is that the Birkhoff factorization takes on
a special form in this case. Recall that Birkhoff factorization applied to the fundamental solution
$S$ (see the appendix) produces
\begin{equation}
\label{bf}
S(\hbar,\hbar^{-1})=Q^{-1}(\hbar)R(\hbar^{-1}).
\end{equation}
From what we mentioned above, in our case the fundamental solution can be taken as
\ba
S^t=
\begin{pmatrix}
I|_{p=0}& \frac{d}{dp}I|_{p=0}
\end{pmatrix}
\ea
Now, the Birkhoff factorization (\ref{bf}) can be rewritten more precisely as
\ba
S(q,\hbar,\hbar^{-1})=(I_2+\sum_{k >0}Q_k q^k)^{-1}R(q,\hbar^{-1})
\ea
where $I_2$ is the 2 dimensional identity matrix and we have set
\ba
A=\sum_{k \ge 0}A_k q^k
\ea
for a matrix $A$.
Looking back at the $I$ function $I_1$ of Section 3.1, it is
evident that in the expansion $\nu=\infty$, the only term contributing
positive powers of $\hbar$ to $I_1$ is 
\ba
\frac{1}{\prod_{m=1}^d (p+\nu+m\hbar)}.
\ea
This implies that the $Q$ matrix of the Birkhoff factorization for this
case takes on the following special form:
\begin{equation}
\label{qfact}
Q_k=Q_k\big((\hbar+\nu)^{-1},(2\hbar+\nu)^{-1} \dots (k\hbar+\nu)^{-1}\big).
\end{equation}
In other words, $Q$ makes sense when expanded as a power series in
either $\hbar$ or $1/\hbar$!

We can use this to immediately write down a function which posesses
$\Xi_3=K_{\p^2}$ and $\Xi_1=\oo(1)\oplus \oo(-3)\rightarrow \p^1$
invariants. (Note that this can readily be modified to include $\Xi_2$
invariants as well, via the comment at the end of section 3.2.) This is
more than a bit surprising, as $\Xi_1$ requires Birkhoff factorization while $\Xi_3$ does
not, i.e. the asymptotic expansion of $\nu$ for each respective space
appears incompatible. At any rate, noting our well-behaved $Q$ above, we simply modify the factorization by inserting a
parameter $w$:
\ba
S'(q,\hbar,\hbar^{-1},w)=(I_2+w\sum_{k >0}Q_k q^k)^{-1}R(q,\hbar^{-1}).
\ea
We now treat $S'$ as a function of $1/\hbar$ and expand. The result is a
new $I$ function 
\ba
I'=1+\frac{pt+t_0}{\hbar}+\frac{W}{\hbar^2}+\dots
\ea
where everything is now $w$ dependent. The mirror maps are
\ba
t_0&=&-2q_1+(17-2w)q_1^2+\Big(50w-\frac{710}{3}\Big)q_1^3+\Big(\frac{8049}{2}-1137w\Big)q_1^4+\dots\\
t&=&(2w-8)q_1+(74-29w)q_1^2+\Big(\frac{1532}{3}w-\frac{3212}{3}\Big)q_1^3+\Big(18609-\frac{19893}{2}w\Big)q_1^4+\dots
\ea
from which we can easily make a consistency check, in that $w=0$ gives
the mirror map for $\Xi_1$ while $w=1$ is the mirror map for
$K_{\p^2}$. After inverting this and changing variables, we find 
\ba
W|_{p=0}=x_1-\Big(2w+\frac{7}{4}\Big)x_1^2+\Big(6w^2+15w+\frac{55}{9}\Big)x_1^3-\Big(\frac{64}{3}w^3+100w^2+\frac{643}{6}w+\frac{455}{16}\Big)x_1^4+\dots
\ea
where we have set $p=0$ to avoid dealing with the leading K\"ahler
parameter factor\footnote{What happens is
that for $w=0$, we get a prefactor of
$(3p+\nu)(p+\nu)|_{p^2=0}=4p\nu+\nu^2$, while $w=1$ gives
$(3p+\nu)^2|_{p^2=0}=6p\nu+\nu^2$. This simply reflects the difference
in the bundle structure of the two spaces.}. Then one sees immediately that $w=0$ gives $\Xi_1$ invariants while $w=1$
gives $K_{\p^2}$ invariants (up to an overall multiplicative factor of $1/3$). By combining this result with those of the
previous two sections, we arrive at a single generating function for the
Gromov-Witten invariants of all $\Xi_i$.

\subsection{A second approach to $\Xi_2$}

We provide here a second derivation of mirror symmetry for $\Xi_2$,
which is along the lines of \cite{FJ4}. This calculation allows us to
provide a closed form for the prepotential and mirror map of $\Xi_2$,
similar to what was found in \cite{FJ4}. Note that we apply the
cohomology relation $p^2=0$ to all toric vectors in this section.

The starting point is the (general) equivariant mirror
symmetry calculation of $\oo(1)\oplus \oo(-3)\rightarrow \p^1$ with
equivariant weights $(1,-\nu)$ acting on the bundle factors, which was
carried out in \cite{FJ3}. Note that these  are the most general equivariant
Gromov-Witten invariants that can be computed for this space. This can be represented by charge vectors
\ba
\label{geninv}
\begin{pmatrix}
1&1&1&-3\\
&&1&-\nu
\end{pmatrix}
\ea
(Note that the 1 in the second row should strictly speaking be set to
$\mu$, an equivariant parameter, but we have set this $\mu=1$ since the
invariants computed using (\ref{geninv}) are already the most general.)
By mirror symmetry, we can use this to compute generalized equivariant
Gromov-Witten invariants, and the potential turns out to be
\ba
\label{gsp}
W(\nu,x_1)=x_1-\Big(\frac{3}{4}\nu+\frac{3}{4}\nu^2+\frac{1}{4}\nu^3\Big)x_1^2+\Big(\frac{13}{12}\nu^2+\frac{9}{4}\nu^3+\frac{23}{12}\nu^4+\frac{3}{4}\nu^5+\frac{1}{9}\nu^6\Big)x_1^3+\dots
\ea
Then $\nu=1$ gives physical GW invariants (i.e. the ones from the
earlier section), while $\nu=-1$ gives the multiple cover formula (see \cite{FJ3}).

Now, one might initially hope that a clever choice of $\nu$ would allow
us to match this to the 2 curve prepotential (up to an overall factor
of 2). Solving the obvious equation yields a value $\nu=12^{1/3}-1$, and
after re-substituting this into (\ref{gsp}), the result does turn out to
be very, very close to the right answer. Unfortunately, there is an error term
present that grows with increasing degree, and without the exact
form of the error, the extraction of integer invariants is not
possible. Nonetheless, this gives us a hint as to how to proceed: one
might very well hope that some
refinement of the theory (\ref{geninv})
contains the looked-for invariants.

A quick glance at \cite{FJ4} reveals that there is indeed a perfectly
natural refinement of (\ref{geninv}), along the lines of the previous section. From \cite{FJ4}, it was shown that the above $\nu$
dependent invariants can also be computed by using  the charge vectors
for the $D_4$ singularity
\ba
\label{d4flat}
\begin{pmatrix}
1&1&1&-1&-2\\
&&1&-\nu&-\nu
\end{pmatrix},
\ea
which holds since the equivariant Euler classes of the bundles are the
same, $-3p-\nu=(-2p-\nu)(-p-\nu)|_{p^2=0}$.
Then the most natural refinement of the equivariant theory defined by (\ref{d4flat}) is of course 
\ba
\label{type1}
\begin{pmatrix}
1&1&1&-1&-2\\
&&1&-\nu_1&-\nu_2
\end{pmatrix}
\ea
Then miraculously, with the simple specialization of weights
$\nu_1=2, \nu_2=1$, we immediately recover the expansion for the
$W_{\Xi_2}$ function!

Notice that this space also naturally contains mirror symmetry for
$\Xi_1=\oo(1)\oplus \oo(-3)\rightarrow \p^1$. This means we can simultaneously carry out mirror symmetry for $\Xi_1$ and
$\Xi_2$. Namely, mirror symmetry for $\Xi_1$ is given by (here we have
restored $\mu$, which was set to 1 in the above calculation)
\ba
\begin{pmatrix}
1&1&1&-1&-2\\
&&\mu&-\mu&-\mu
\end{pmatrix},
\ea
and hence mirror symmetry for $\Xi_1$ and $\Xi_2$ is determined by 
\ba
\begin{pmatrix}
1&1&1&-1&-2\\
&&\mu&-\mu-\nu&-\mu
\end{pmatrix}.
\ea
That is, $\nu=0$ gives $\Xi_1$ while $\nu=\mu$ is $\Xi_2$.

We note briefly that the $\Xi_2$ invariants can also be computed from the
toric data 
\ba
\label{second}
\begin{pmatrix}
1&1&1&-1&-1&-1\\
&&1&-2\nu&-\nu&-\nu
\end{pmatrix},
\ea
which again follows since the equivariant Euler classes of both cases are the
same. We will come back to this later in the paper.

\section{The localization calculation}

In this section, we reproduce the refined Gromov-Witten invariants of $K_{\p^2}$ derived in 
\cite{AMV} by using the standard localization computation with specialized torus action weights on $\p^2$.
The localization computation starts by assigning torus action
weights to the homogeneous coordinates of $\p^2$:
\begin{equation}
(X_{1}:X_{2}:X_{3})\mapsto (e^{\lambda_{1}t}X_{1}:e^{\lambda_{2}t}X_{2}:e^{\lambda_{3}t}X_{3}).
\end{equation}
In the standard computation of local Gromov-Witten invariants of $K_{\p^2}$, we obtain the 
same results for an arbitrary choice of weights
$(\lambda_{1},\lambda_{2},\lambda_{3})$, but to obtain the
refined invariants of \cite{AMV}, we specialize these weights as follows:
\begin{equation}
(\lambda_{1}, \lambda_{2},\lambda_{3})=(1,\omega,\omega^{2}),
\label{basic}
\end{equation}
where $\omega$ is the primitive cubic root of unity $-\frac{1}{2}+\frac{\sqrt{-3}}{2}$.
Of course, we can permute the subscripts of weights, and the following operations is also allowed:
\begin{eqnarray}
&& (1, \omega,\omega^{2}) \mapsto (e^{\lambda}, e^{\lambda}\omega,e^{\lambda}\omega^{2}),
\end{eqnarray}
because these operations keep each amplitude $a(\Gamma)$, associated
with the tree graph $\Gamma$ used in the localization
computation of invariants, invariant. But any choice of weights which
changes $a(\Gamma)$ is no longer allowed when computing 
refined invariants in \cite{AMV}. The reason for this restriction comes from the following 
operation. The tree graph $\Gamma$ used in the localization computation of $K_{\p^2}$ consists of 
colored vertices $v(1), v(2), v(3)$ associated with homogeneous coordinates $X_{1},X_{2},X_{3}$ and of edges 
with degree $d$. Then the degree $d$ genus $0$ local Gromov-Witten invariant of $K_{\p^2}$ is given
by the following formula:
\begin{equation}
\sum_{{\Gamma},\;{d(\Gamma)=d}}a(\Gamma)q^{d(\Gamma)},
\label{sum}
\end{equation}
where $q$ is the degree counting parameter and $d(\Gamma)$ is the degree of the 
graph $\Gamma$ defined by the sum of degrees of edges in $\Gamma$. But to obtain the refined invariants 
of $K_{\p^2}$ in \cite{AMV}, we have to introduce three different degree counting parameters
$x_{1},x_{2},x_{3}$ that correspond to the three three $\p^1$'s in the toric diagram of $K_{\p^2}$. 
To realize this refinement, what we have to do is to define a colored degree $(d_{1}(\Gamma),d_{2}(\Gamma),d_{3}(\Gamma))$
for the graph $\Gamma$. $d_{1}(\Gamma)$ (resp. $d_{2}(\Gamma)$, $d_{3}(\Gamma)$) is defined as the sum of 
degrees of edges whose boundary vertices consist of $v(2)$ and $v(3)$ (resp. $v(1)$ and $v(3)$, $v(1)$ and $v(2)$).
Then we change (\ref{sum}) as follows:
\begin{equation}
\sum_{{\Gamma},\;{d(\Gamma)=d}}a(\Gamma)x_{1}^{d_{1}(\Gamma)}x_{2}^{d_{2}(\Gamma)}x_{3}^{d_{3}(\Gamma)}.
\label{rsum}
\end{equation}
If we set the weights as proposed in (\ref{basic}), the above formula gives the generating function of refined 
invariants of \cite{AMV} of total degree $d$. But if we set the weights
in a way not allowed according to previous discussion, 
(\ref{rsum}) turns out to be a meaningless polynomial. In the following,
we write down the results of the computation 
in the form of a generating function up to total degree $6$:  
\begin{eqnarray}
F(w)&:=& \biggl({x_{3}} + {x_{1}} + {x_{2}}\biggr)w \no\\
&&+ \biggl( - {x_{1}}{x
_{2}} - { \frac {7}{8}} {x_{3}}^{2} - {x_{2}}{x
_{3}} - { \frac {7}{8}} {x_{2}}^{2} - {x_{3}}{x
_{1}} - { \frac {7}{8}} {x_{1}}^{2}\biggr)w^2 \no\\
&&+ \biggl(3{x_{1}}{x_{3}}^{2} + 3{x_{2}}^{2}{x_{3}}
 + 3{x_{2}}{x_{3}}^{2} + 3{x_{1}}{x_{2}}^{2} + 
{ \frac {55}{27}} {x_{1}}^{3} + { 
\frac {55}{27}} {x_{3}}^{3} + 3{x_{1}}{x_{2}}{x_{3}} + 3
{x_{1}}^{2}{x_{3}} \no\\
&&+{ \frac {55}{27}} {x_{2}}
^{3} + 3{x_{1}}^{2}{x_{2}}\biggr)w^3  \no \\
&&+\biggl( - 13{x_{1}}{x_{2}}^{3} - 13{x_{1}}{x_{3}}^{3} - 16{x
_{1}}{x_{2}}{x_{3}}^{2} - { \frac {121}{8}} {x
_{1}}^{2}{x_{3}}^{2} - { \frac {455}{64}} {x_{2}
}^{4} - { \frac {121}{8}} {x_{1}}^{2}{x_{2}}^{2}\no
 \\
&& - 13{x_{1}}^{3}{x_{2}} - 13{x_{2}}^{3}{x_{3}} - 
16{x_{1}}^{2}{x_{2}}{x_{3}} - 16{x_{1}}{x_{2}}^{2}{x
_{3}} - 13{x_{1}}^{3}{x_{3}} - { \frac {455}{64}
} {x_{3}}^{4} \no\\
&& - { \frac {121}{8}} {x_{2}}^{2}{x_{3}}^{
2} - { \frac {455}{64}} {x_{1}}^{4} - 13{x_{2}}
{x_{3}}^{3}\biggr)w^4 \no\\
&&+ \biggl(104{x_{1}}{x_{2}}{x_{3}}^{3}
 + 68{x_{1}}^{4}{x_{2}} + 91{x_{2}}^{3}{x_{3}}^{2}
 + { \frac {3876}{125}} {x_{2}}^{5} + 91{
x_{2}}^{2}{x_{3}}^{3} + 68{x_{2}}^{4}{x_{3}} \no\\
&&+ 112{x_{1}}
{x_{2}}^{2}{x_{3}}^{2} + 68{x_{1}}{x_{3}}^{4} + 112{x_{
1}}^{2}{x_{2}}{x_{3}}^{2}+ 91{x_{1}}^{2}{x_{3}}^{3} + 68{x_{1}}^{4}{x_{3}}
 + 91{x_{1}}^{3}{x_{2}}^{2} \no\\
&&+ 112{x_{1}}^{2}{x_{2}}^{2}
{x_{3}} + 91{x_{1}}^{2}{x_{2}}^{3} + 68{x_{2}}{x_{3}}^{4}
 + 91{x_{1}}^{3}{x_{3}}^{2} + 104{x_{1}}{x_{2}}^{3
}{x_{3}} + 104{x_{1}}^{3}{x_{2}}{x_{3}}\no\\ 
&&+ 68{x_{1}}{x
_{2}}^{4} + { \frac {3876}{125}} {x_{1}}^{5} + 
{ \frac {3876}{125}} {x_{3}}^{5}\biggr)w^5 \no\\
&&+\biggl(- 399{x_{2}}^{5}{x_{3}} - { \frac {4845}{8}} 
{x_{1}}^{2}{x_{3}}^{4} - 891{x_{1}}^{2}{x_{2}}{x_{3}}^{3}
 - 891{x_{1}}{x_{2}}^{3}{x_{3}}^{2} - 891{x_{1}}{x_{2}}
^{2}{x_{3}}^{3}\no \\
&& - 741{x_{1}}{x_{2}}{x_{3}}^{4} - { 
\frac {18496}{27}} {x_{1}}^{3}{x_{3}}^{3} - { 
\frac {4845}{8}} {x_{1}}^{4}{x_{2}}^{2} - { 
\frac {18496}{27}} {x_{1}}^{3}{x_{2}}^{3} - { 
\frac {4845}{8}} {x_{1}}^{2}{x_{2}}^{4}\no \\
&& - 399{x_{1}}{x_{2}}^{5} - 399{x_{1}}^{5}{x_{2}}
 - 741{x_{1}}^{4}{x_{2}}{x_{3}} - 891{x_{1}}^{3}{x_{2}}
{x_{3}}^{2} - { \frac {33649}{216}} {x_{2}}^{6}
 - { \frac {33649}{216}} {x_{3}}^{6}\no \\
&& - { \frac {33649}{216}} {x_{1}}^{6} - 891
{x_{1}}^{2}{x_{2}}^{3}{x_{3}} - { \frac {7533
}{8}} {x_{1}}^{2}{x_{2}}^{2}{x_{3}}^{2} - 891{x_{1}}^{3}
{x_{2}}^{2}{x_{3}} - { \frac {4845}{8}} {x_{2}
}^{4}{x_{3}}^{2}\no \\
&& - { \frac {4845}{8}} {x_{2}}^{2}{x_{3}}
^{4} - 399{x_{2}}{x_{3}}^{5} - 399{x_{1}}^{5}{x_{3}} - 
399{x_{1}}{x_{3}}^{5} - { \frac {4845}{8}} {x
_{1}}^{4}{x_{3}}^{2} - { \frac {18496}{27}} {x_{
2}}^{3}{x_{3}}^{3}\no \\
&& - 741{x_{1}}{x_{2}}^{4}{x_{3}}\biggr)w^6+\cdots. 
\end{eqnarray}
Note that the notation here differs from previous sections; above, $x_i$
are the equivariant parameters corresponding to the $\mu_i$ of Section
2, and $w$ is the same as the $x_1$ of Section 3.
If we set $x_{3}=0$, we obtain generating function of local Gromov-Witten invariants 
of the two curves case with two K\"ahler parameters:
\begin{eqnarray}
F|_{x_{3}=0}(w) &:= &\biggl({x_{2}} + {x_{1}}\biggr)w 
+ \biggl( - { 
\frac {7}{8}} {x_{2}}^{2} - {x_{1}}{x_{2}} - { 
\frac {7}{8}} {x_{1}}^{2}\biggr)w^2+ \biggl({ \frac {55}{
27}} {x_{1}}^{3} + { \frac {55}{27}} {x_{2}}^{3}
 + 3{x_{1}}{x_{2}}^{2} + 3{x_{1}}^{2}{x_{2}}\biggr)w^3\no \\
&& + \biggl( - 13{x_{1}}{x_{2}}^{3} - { \frac {
455}{64}} {x_{1}}^{4} - 13{x_{1}}^{3}{x_{2}} - 
{ \frac {121}{8}} {x_{1}}^{2}{x_{2}}^{2} - 
{ \frac {455}{64}} {x_{2}}^{4}\biggr)w^4 \no\\
&& + \biggl(68{x_{1}}{x_{2}}^{4} + { \frac {3876
}{125}} {x_{1}}^{5} + 91{x_{1}}^{3}{x_{2}}^{2} + 
{ \frac {3876}{125}} {x_{2}}^{5} + 68{x_{1}}^{4}
{x_{2}} + 91{x_{1}}^{2}{x_{2}}^{3}\biggr)w^5 \no\\
&&+ \biggl( - 399{x_{1}}{x_{2}}^{5} - { \frac {18496}{27}} 
{x_{1}}^{3}{x_{2}}^{3} - { \frac {4845}{8}} {x
_{1}}^{4}{x_{2}}^{2} - { \frac {33649}{216}} {x
_{1}}^{6} - { \frac {33649}{216}} {x_{2}}^{6} - 
{ \frac {4845}{8}} {x_{1}}^{2}{x_{2}}^{4}\no \\
&& - 399{x_{1}}^{5}{x_{2}}\biggr)w^6+\cdots .
\end{eqnarray}
If we set $x_{1}=x_{2}=1,x_{3}=0$, we obtain one parameter generating function of the two curves case: 
\begin{eqnarray}
F|_{x_{1}=x{_2}=1,x_{3}=0}(w):= 2w - { \frac {11}{4}} w^{2} + 
{ \frac {272}{27}} w^{3} - { \frac {
1771}{32}} w^{4} + { \frac {47502}{125}} w^{5}
 - { \frac {81158}{27}} w^{6}+\cdots.
\end{eqnarray}
We can easily see that $w\frac{d}{dw}F|_{x_{1}=x_{2}=1,x_{3}=0}$ agrees with the result 
from the mirror symmetry computation.
Moreover, if we set $x_{2}=x_{3}=0$, the result agrees with the one ${\cal O}(1)\oplus{\cal O}(-3)$ curve:  
\begin{eqnarray}
F|_{x_{1}=1,x_{2}=x_{3}=0}(w) := {w} - { \frac {7}{8}} {w}^{2} + { \frac {55}{27}} {w}^{3} - 
{ \frac {455}{64}} {w}^{4} + 
{ \frac {3876}{125}} {w}^{5} - 
{ \frac {33649}{216}} {w}^{6}+\cdots.
\end{eqnarray}

\section{Revisiting and generalizing $\oo(k)\oplus
 \oo(-2-k)\rightarrow \p^1$}

\subsection{A new way to compute mirror symmetry for $\oo(k)\oplus
 \oo(-2-k)\rightarrow \p^1$}

The second main focus of this paper is to present an alternative formulation for the notion of mirror symmetry on
$X_k=\oo(k)\oplus \oo(-2-k)\rightarrow \p^1$ proposed in
\cite{FJ4}. This new viewpoint allows us to give a closed form for the
mirror map for all $k$, and moreover produces the closed form for GW
invariants found in \cite{M}. This approach is closer to physical mirror
symmetry than the formulation
of \cite{FJ4}, in the sense that we can use a single $I$ function for
the computation for all $k$.

We first recall the work of \cite{FJ4}. The main observation of that
paper was that the Gromov-Witten theory of $X_k$ with antidiagonal
equivariant action $(\mu,-\mu)$ on the bundle is equivalent to GW theory
on the space
\ba
\label{bigbundle}
\oo(1)^{\oplus k}\oplus \oo(-1)^{\oplus (2+k)}\rightarrow \p^1
\ea
with torus action $(\underbrace{\mu \dots \mu}_k, \underbrace{-\mu \dots -\mu}_{2+k})$ on the
bundle. The essential thing we needed to notice to see that these theories
should be the same is that the equivariant Euler classes the respective bundles agree:
\ba
(kp+\mu)((-2-k)p-\mu)|_{p^2=0}=(p+\mu)^k(-p-\mu)^{2+k}|_{p^2=0}.
\ea

So, for example, consider $X_2=\oo(2)\oplus \oo(-4)\rightarrow \p^1$
with action $(\mu,-\mu)$ on the bundle. This can be represented
torically as
\ba
\bmm
1&1&2&-4\\
&&\mu&-\mu
\emm
\ea
where the second row gives the equivariant weights acting on the bundle
factors represented by $2,-4$. Then the results of \cite{FJ4} show that
the same Gromov-Witten theory can be calculated by using the toric data
\ba
\begin{pmatrix}
1&1&1&1&-1&-1&-1&-1\\
&&\mu&\mu&-\mu&-\mu&-\mu&-\mu
\end{pmatrix}
\ea
which corresponds to the space
\ba
\oo(1)_{\mu}\oplus\oo(1)_{\mu}\oplus\oo(-1)_{-\mu}\oplus\oo(-1)_{-\mu}\oplus\oo(-1)_{-\mu}\oplus\oo(-1)_{-\mu}\rightarrow \p^1.
\ea
As usual, subscripts denote equivariant weights on the respective
factors. The point of doing this reduction in \cite{FJ4} is that this allows us
to find a very simple mirror map which agrees with \cite{M}, as well as
to find a closed form for the Yukawa coupling.

Naturally, since only the agreement of equivariant Euler classes is
needed for the above theories to coincide, the choice of bundle
$\oo(1)^{\oplus 2}\oplus \oo(-1)^{\oplus 4}$ is not unique,
so there may be any number of (possibly even simpler) theories which
share the same GW invariants. One finds immediately that the
equivariant Euler class of the bundle on $X_k$ is also the same on the space
\ba
\label{newcurve}
\begin{pmatrix}
1&1&1&-1&-1&-1\\
&&\mu/k&-\mu/k&-\mu&-\mu
\end{pmatrix}
\ea
To get some geometric insight into what this vector means, we note that the toric data of equivariant
$\oo(-1)\oplus\oo(-1)\rightarrow \p^1$ with weights $(-\mu,-\mu)$ acting
on the bundle is
\ba
\begin{pmatrix}
1&1&-1&-1\\
&&-\mu&-\mu
\end{pmatrix}
\ea
and the difference between this and (\ref{newcurve}) is just the
middle  two columns, which represent the bundles $\oo(1)\oplus
\oo(-1)$. Thus, this can be thought of as a `twist' in the sense of
Coates-Givental \cite{CG}, which geometrically means that
(\ref{newcurve}) is the total space
\ba
\oo(1)_{\mu/k}\oplus\oo(-1)_{-\mu/k}\longrightarrow \Big(\oo(-1)_{-\mu}\oplus \oo(-1)_{-\mu}\rightarrow \p^1\Big).
\ea
Hence, $\oo(k)\oplus \oo(-2-k)\rightarrow \p^1$ for all $k$ is just a
twist of the $\oo(-1)\oplus \oo(-1)$ curve.

Call the space defined by (\ref{newcurve})
$X_k'$. Then, the advantage of working with $X_k'$ rather than
(\ref{bigbundle}) is clear; the fact that the rank of the bundle of $X_k'$ doesn't
change (i.e., the bundle of (\ref{newcurve}) is rank 4 for all $k$,
whereas (\ref{bigbundle}) has rank $2+2k$ ) allows us to use a single $I$ function to compute the answer for
all $k$ (in contrast to (\ref{bigbundle})). To be completely explicit,
the $I$ function  corresponding to (\ref{bigbundle}) for arbitrary $k$
is given as
\ba
\label{i1}
I_1=q^{p/\hbar}\sum_{d \ge 0}\frac{\prod_{m=-d+1}^0(-p-\mu+m\hbar)^{2+k}}{\prod_{m=1}^d(p+\mu+m\hbar)^{k}\prod_{m=1}^d(p+m\hbar)}q^d
\ea
while that for (\ref{newcurve}) is
\ba
\label{i2}
I_2=q^{p/\hbar}\sum_{d \ge 0}\frac{\prod_{m=-d+1}^0(-p-\mu/k+m\hbar)\prod_{m=-d+1}^0(-p-\mu+m\hbar)^2}{\prod_{m=1}^d(p+\mu/k+m\hbar)\prod_{m=1}^d(p+m\hbar)}q^d.
\ea
Then it is evident that $I_2$ can be Birkhoff factorized for arbitrary
values of $k$, whereas $I_1$ cannot; thus we can use $I_2$ to find the
general form of the mirror transformation for $X_k$.

Then, by going through the usual routine, from the above vector
(\ref{newcurve}) we can expand in $\mu=\infty$ and perform Birkhoff
factorization in order to extract the $J$ function for arbitrary values
of $k$:
\ba
J_k=q^{p/\hbar}\Big(1+\frac{\mu
k\big(p(2+k)+1\big)\log(1+q)}{\hbar}+\frac{\tilde W_k(q)}{\hbar^2}+\dots
\Big)
\ea
i.e. the mirror map is 
\ba
t_k(q)=\log(q)+k(2+k)\log(1+q),
\ea
in agreement with \cite{FJ4} and \cite{M}. Then, after inverting the
mirror map, we find
\ba
\tilde W_k(q(t))=\mu(2pk+2p+\mu)W_k
\ea
where remarkably, the coefficients of $W_k$ can easily be seen to have a
very simple \it closed form \normalfont:
\ba
W_k=\sum_{j=1}^{\infty}\frac{(-1)^{(j-1)k}}{j^2(j-1)!}\prod_{m=1}^{j-1}(k(2+k)j+m)x^k
\ea
where $x=e^t$. This is, of course, exactly the prepotential of \cite{M},
and, as we will see in the next section, a slight modification of this
will allow us to deduce a closed form for the invariants of $\Xi_2$.

\subsection{More general equivariant theories}

We now turn to a generalization of the above equivariant theory which,
 surprisingly, also turns out to have an elegant closed form, and
 further contains $\Xi_2$ invariants, which is how we derived the
 closed formula presented earlier.

The idea here is the
following. 
From \cite{FJ3}, mirror symmetry for $X_1$ can be performed by use of
the vector
\ba
\begin{pmatrix}
1&1&1&-3\\
&&\mu&-\mu
\end{pmatrix}.
\ea
Then, from \cite{FJ4}, this can also be computed with
\ba
\label{3spread}
\begin{pmatrix}
1&1&1&-1&-1&-1 \\
&&\mu&-\mu&-\mu&-\mu
\end{pmatrix}
\ea

Now, once we have uncovered the toric model (\ref{3spread}), it is of
course very natural to ask whether there might be some meaning to the
equivariant Gromov-Witten theory of a generalized setup:
\ba
\label{gmod}
\begin{pmatrix}
1&1&1&-1&-1&-1 \\
&&\nu_1&-\nu_2&-\nu_3&-\nu_4
\end{pmatrix}
\ea
Computationally speaking, this model simplifies dramatically if we specialize
the weights to
\ba
\label{gmod2}
\begin{pmatrix}
1&1&1&-1&-1&-1 \\
&&\mu&-\mu&-\nu_1&-\nu_2
\end{pmatrix}
\ea
The reason for this is, the third and fourth columns of (\ref{gmod})
contribute the following factors to the $I$ function:
\ba
\frac{\prod_{m=-d+1}^0(-p-\nu_2+m\hbar)}{\prod_{m=1}^d(p+\nu_1+m\hbar)}.
\ea
Then, almost all the terms on the top and bottom will cancel if we set
$\nu_1=\nu_2$ (that is,
exactly one term on the top and one term on the bottom will
survive). Hence, we consider instead the theories defined by (\ref{gmod2}).

To simplify the result, we take $\mu=1$. Let $J_{\nu}$ be the $J$
function obtained by Birkhoff factorization of (\ref{gmod2}). Amazingly,
the result is almost the same as the previous section:
\ba
J_{\nu}=q^{p/\hbar}\Big(1+\frac{\big(p(\nu_1\nu_2+\nu_1+\nu_2)+\nu_1\nu_2\big)\log(1+q)}{\hbar}+\frac{\tilde W_{\nu}(q)}{\hbar^2}+\dots
\Big)
\ea
Moreover, we once again find a closed form for the Gromov-Witten invariants:
\ba
W_{\nu}(x)=\sum_{k>0}\frac{(-1)^{k+1}}{k^2(k-1)!}\prod_{j=1}^{k-1}(k(\nu_1\nu_2+\nu_1+\nu_2)+j)x^k.
\ea
Here
\ba
\tilde W_{\nu}(q(t))=\big(p(2\nu_1\nu_2+\nu_1+\nu_2)+\nu_1+\nu_2 \big)W_{\nu}(q(t)).
\ea
Now, from (\ref{gmod2}) we clearly see that
\ba
W_{(\nu_1=1,\nu_2=1)}=W_{\Xi_1},
\ea
i.e. the invariants  agree with $\oo(1)\oplus \oo(-3)\rightarrow \p^1$,
as they must. Also, as we saw above,
\ba
W_{(\nu_1=1,\nu_2=2)}=W_{\Xi_2},
\ea
which provides us with a closed formula for the Gromov-Witten invariants
of $\Xi_2$.

\subsection{Connection  open strings on $\C^3$}

\begin{figure}[t]
\label{onelag}
\centering
\input{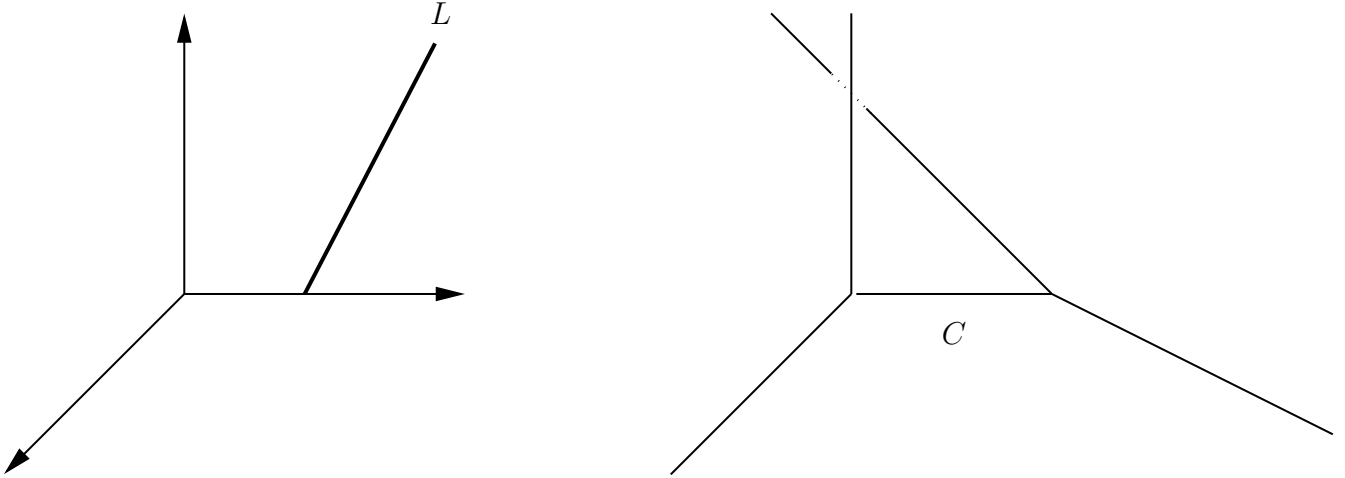}
\caption{The comparison between $\C^3$ with one Lagrangian submanifold
 and a $\p^1$ with normal bundle $\oo(k)\oplus \oo(-2-k)$. Since the
 angle $L$ makes with the base of $\C^3$ is determined by the framing,
 we can think of the framing as determining the bundle structure of the
 right hand side picture. Note that the position of $L$ on the toric
 diagram of $\C^3$ maps to the K\"ahler parameter on $C$.}
\end{figure}

\begin{figure}[t]
\label{2lags}
\centering
\input{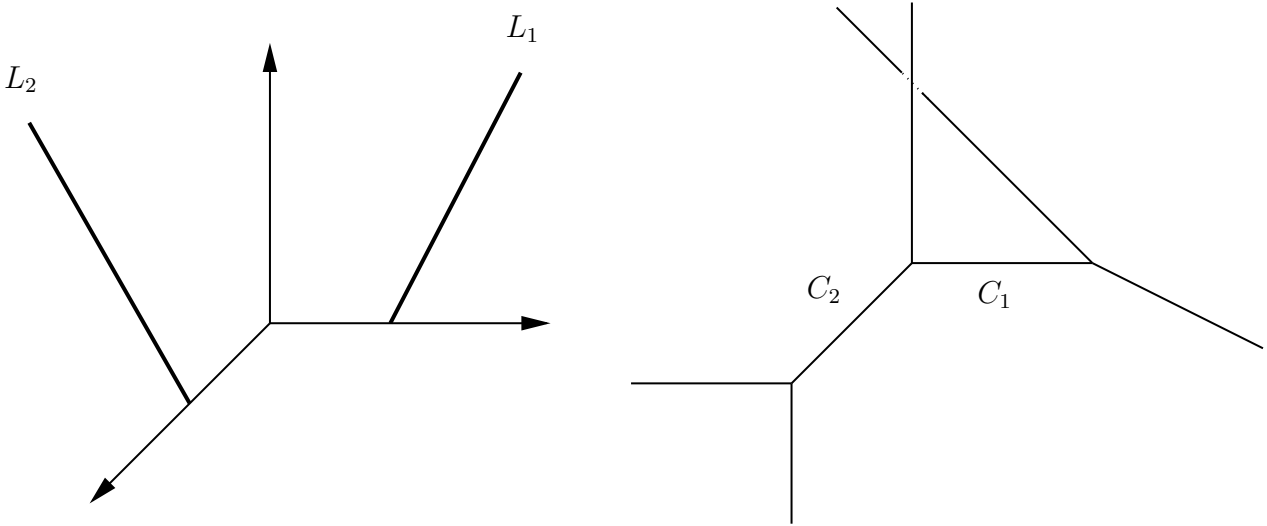}
\caption{The comparison between two Lagrangian submanifolds in $\C^3$
 and the case of two adjacent curves with arbitrary normal bundles (here, we have
 drawn
 $(-1,-1)$ and $(1,-3)$). Again, one can see a formal similarity between
 the two cases (for appropriate choices of framing on $L_1, L_2$).}
\end{figure}

We now outline the connection between closed strings  on
$\oo(k)\oplus\oo(-2-k)\rightarrow \p^1$ and open strings on
$\C^3$. Besides being of independent interest, we hope that eventually a generalization of this correspondence will allow us to define
a two K\"ahler parameter toric model for $\Xi_2$.

We briefly recall the toric approach to open string theory given in
\cite{LMW}. Let $X$ be a noncompact Calabi-Yau threefold, which without
loss of generality we take
to have one K\"ahler modulus. As usual, we represent $X$ as a toric
quotient
\ba
X=\{ z \in \C^4: \sum_i l_i|z_i|^2=r \}/S^1
\ea
where
\ba
S^1: z_i \longrightarrow e^{l_i \theta} z_i.
\ea
Clearly, the vector $l=(l_i)$ completely specifies $X$. Then, according
to \cite{LMW}, open strings on $X$ with arbitrary framing $n$ can be
represented by the vectors
\ba
\begin{pmatrix}
l_1&l_2&l_3&l_4&0&0\\
1&-n-1& n&0&1&-1
\end{pmatrix}
\ea
which is of course the toric data of a CY fourfold. The idea behind this
construction is that the second vector corresponds to the position of a
Lagrangian submanifold at some point on the skeleton of the toric
diagram of $X$, and this subsequently has an interpretation as a cycle
in M theory.

Now, if one were to consider open strings on $\C^3$, since $\C^3$ has no
K\"ahler modulus, the toric data of a single Lagrangian submanifold in
$\C^3$ will be just
\ba
\bmm
1&-n-1&n&1&-1
\emm
\ea
The point is that if one looks at the toric diagram of this open string
geometry (Figure 3), it looks almost the same as the diagram of a $\oo(k)\oplus
\oo(-2-k)\rightarrow \p^1$ curve; this is the correspondence that was
first observed at trivial framing in \cite{AV}, and subsequently
generalized to arbitrary framing in \cite{M}. Then, it is a simple matter to
use \cite{FJ4} to show that these theories are in fact exactly the
same for arbitrary framing.  All we have to do is consider the equivariant version of the
single-Lagrangian-in-$\C^3$ theory:
\ba
\label{openc3}
\bmm
1&1&n&-n-1&-1\\
&&\mu&-\mu&-\mu
\emm
\ea
We then just have to compare this to the equivariant theory on $\oo(k)\oplus
\oo(-2-k)\rightarrow \p^1$:
\ba
\label{curve}
\begin{pmatrix}
1&1&k&-2-k\\
&&\mu&-\mu
\end{pmatrix}
\ea
Then, as was shown in \cite{FJ4}, the data of (\ref{curve}) is
equivalent to that of 
\ba
\label{curve2}
\bmm
1&1&k&-k-1&-1\\
&&\mu&-\mu&-\mu
\emm,
\ea
which is true simply because the equivariant Euler classes of the
bundles of
(\ref{curve}) and (\ref{curve2}) are the same.

Then clearly, (\ref{openc3}) and (\ref{curve2}) agree, which shows
quantitatively the correspondence between open strings on $\C^3$ with
arbitrary framing and closed strings on $\oo(k)\oplus
\oo(-2-k)\rightarrow \p^1$.

Finally, we mention that the 2 curve problem $\Xi_2$, which consists of 2
curves with normal bundle $\oo(1)\oplus \oo(-3)$ intersecting at right
angles, is supposed to be equivalent to the GW theory of open strings on
$\C^3$ with two Lagrangian submanifolds on adjacent legs \cite{M2}. Again, by
looking at the toric diagram, we can see the formal similarity between
these theories. However, it is not clear how one would define a toric
model and Birkhoff factorization (which is required to compute the
invariants in all the examples of this section) for a 2 parameter case.

\section{Conclusion}

In this paper, we have taken the first steps toward establishing a fully
equivariant version of mirror symmetry for $K_{\p^2}$, by considering
the special case a single independent equivariant parameter. We have also
given a novel approach to localization which reproduces the fully
equivariant invariants of \cite{AMV}. This represents a major improvement over the
results of \cite{FJ4}, in that $K_{\p^2}$ has a nontrival 4 cycle, which,
as we have seen, compounds the difficulty of the problem immensely.

Many questions are raised by our work here. The first and most obvious
is how one could turn the contents of this paper into a complete theory
of equivariant mirror symmetry of $K_{\p^2}$. The complication here is,
as we have seen above, considering mirror symmetry for `subgraphs' of
the toric diagram for $K_{\p^2}$ requires Birkhoff factorization. At the
same time, $K_{\p^2}$'s mirror symmetry requires that we do \it not
\normalfont Birkhoff factorize. We have reconciled this problem in the
present paper, but how this would generalize to a 3 equivariant
parameter case is unclear.

Another troublesome
issues is how to extend these techniques to a case with two K\"ahler
parameters, such as $K_{F_n}$ for the Hirzebruch surfaces $F_n$, even
for a single independent equivariant parameter. It is
not immediately clear how this would work, since even for $K_{F_1}$,
there is an $\oo(1)\oplus \oo(-3)\rightarrow \p^1$ curve in the geometry
which should be visible in the fully equivariant model. We hope to
address these questions in future work.

\section{Appendix: Period integrals, Birkhoff factorization and mirror symmetry}

Here we collect the relevant background material on working out
equivariant mirror symmetry from toric geometry. For simplicity, all
considerations are restricted to the case of one K\"ahler parameter.

We can define a toric manifold $X$ by starting with a vector
\ba
l=
\begin{pmatrix}
l_1&l_2&\dots&l_n
\end{pmatrix},
\ea
where $l_i >0 \ \  \forall i$,
and then using it to write a symplectic quotient
\ba
X=\{(z_i)\in \C^n: \sum_i l_i|z_i|^2=r\}/S^1
\ea
where the $S^1$ action is given as
\ba
S^1:(z_1\dots z_n)\rightarrow (e^{\sqrt{-1} l_1 \theta}z_1 \dots e^{\sqrt{-1} l_n \theta}z_n).
\ea
We can then consider the equivariant version of $X$ by including the
torus action $X$ induced by that on $\C^n$:
\ba
(\C^*)^n:(z_1\dots z_n)\rightarrow (\mu_1z_1 \dots \mu_n z_n).
\ea
We represent this in this paper as a matrix
\ba
\begin{pmatrix}
l_1& \dots& l_n\\
\mu_1&\dots&\mu_n.
\end{pmatrix}
\ea
This gives rise to the equivariant cohomology ring for $X$:
\ba
H^*_T(X,\C)=\frac{\C[p] \otimes \Q[\mu]}{<\prod_i(l_ip+\mu_i)>}
\ea

From this data, we can construct Givental's equivariant $I$ function, which is
supposed to give the period integrals of the mirror manifold to $X$:
\ba
I=q^{p/\hbar}\sum_{d\ge 0}\frac{1}{\prod_{i=1}^n\prod_{m=1}^{l_i d}(l_ip+\mu_i+m\hbar)}q^d.
\ea
This is a series with coefficients taking values in $H^*_T(X,\C)$.

We can then use \cite{CG} to consider mirror symmetry for the total
space of $X$ with line bundles over it. Without loss of generality we
may restrict to the case $\oo(k_1)\oplus \oo(-k_2)\rightarrow X$ for $k_1,k_2
>0$. Then from \cite{CG}, the $I$ function is modified as follows:
\ba
I^{twist}=q^{p/\hbar}\sum_{d\ge 0}\frac{1}{\prod_{i=1}^n\prod_{m=1}^{l_i
d}(l_ip+\mu_i+m\hbar)}\times \frac{\prod_{m=-k_2d+1}^{0}(-k_2 p-\nu_2+m\hbar)}{\prod_{m=1}^{k_1
d}(k_1p+\nu_1+m\hbar)} q^d
\ea 
and we are supposed to expand this twisted function about
$\nu_1=\nu_2=\infty$.

In order to exhibit mirror symmetry, one must now use Birkhoff
factorization, which goes as follows. Upon expansion of $I^{twist}$ in
$\nu_1=\infty$, positive powers of $\hbar$ will be introduced into
$I^{twist}$. We have to remove these via a procedure called `Birkhoff
factorization'. First, define the fundamental solution:
\ba
S(\hbar,\hbar^{-1})=
M \times
\begin{pmatrix}
I^{twist}|_{p=-\mu_1/l_1} &\dots & I^{twist}|_{p=-\mu_n/l_n}
\end{pmatrix}
\ea
(here $M$ is a change of basis matrix which brings $S$ to the basis of
solutions $(I^{twist}|_{p=0}, \dots , \frac{d^n}{dp^n} I^{twist}|_{p=0})
$).
Then this factors as
\ba
S(\hbar,\hbar^{-1})=Q^{-1}(\hbar)R(\hbar^{-1}).
\ea
Then $R(\hbar^{-1})$ gives the `factored' fundamental solution, from
which we can write down the $J$ function, which takes on the form
\ba
J=1+\frac{t(q)}{\hbar}+\frac{W(q)}{\hbar^2}+\dots
\ea
Then $t(q)$ is the mirror map, and we can recover the Gromov-Witten
invariants of the configuration $\oo(k_1)\oplus \oo(-k_2)\rightarrow X$
from the function
\ba
W(q(t)).
\ea

\end{document}